\documentclass[12pt]{article}
\usepackage{amsfonts,amssymb,subeqnarray}
\def\cA{{\cal A}}                    
                    
\def\cG{{\cal G}}                    
                    \def\cL{{\cal L}}

                    \def\cU{{\cal U}}
          \def\cW{{\cal W}}          
\def\cY{{\cal Y}}          

\def\cotan{\mathop{\rm cotan}\nolimits}

\newcommand{\CC}{{\mathbb C}}
\newcommand{\II}{{\mathbb I}}
\newcommand{\ZZ}{{\mathbb Z}}
\newcommand{\eps}{{\varepsilon}}
\newcommand{\elp}{{{\cal A}_{q,p}(\widehat{sl}(N)_{c})}}
\newcommand{\elpa}[1]{{{\cal A}_{q,p}(\widehat{sl}(#1)_{c})}}

\newcommand{\deysr}{{{\cal D}Y(sl(N))_c}}
\newcommand{\dey}{{{\cal D}Y_{r}(sl(N))_c}}
\newcommand{\deygl}{{{\cal D}Y_{r}(gl(N))_c}}
\newcommand{\uq}[2]{{{\cal U}_{#1}(\widehat{sl}(#2)_{c})}}

\newcommand{\sfrac}[2]{{\textstyle{\frac{#1}{#2}}}}
\newcommand{\half}{{\sfrac{1}{2}}}
\newcommand{\car}[2]{\left[\begin{array}{c}{#1}\\{#2}\end{array}\right]}
\newcommand{\finproof}{{\hfill \rule{5pt}{5pt}}}
\newcommand{\wL}{\widetilde L}

\newtheorem{coro}{Corollary}
\newtheorem{lemm}{Lemma}
\newtheorem{thm}{Theorem}

\begin{document}
\pagestyle{empty}
\vfill
\begin{center}

{\Large \textbf{Deformed Double Yangian Structures}}

\vspace{10mm}

{\large D. Arnaudon$^a$, J. Avan$^b$, L. Frappat$^a$, M. Rossi$^c$}

\vspace{10mm}

\emph{$^a$ Laboratoire d'Annecy-le-Vieux de Physique Th{\'e}orique LAPTH}

\emph{CNRS, URA 1436, associ{\'e}e {\`a} l'Universit{\'e} de Savoie}

\emph{LAPP, BP 110, F-74941 Annecy-le-Vieux Cedex, France}

\vspace{7mm}

\emph{$^b$ LPTHE, CNRS, UMR 7589, Universit{\'e}s Paris VI/VII, France}

\vspace{7mm}

\emph{$^c$ Department of Mathematics, University of Durham \\
South Road, Durham DH1 3LE, UK}

\end{center}

\vfill
\vfill

\begin{abstract}
Scaling limits at $q \rightarrow 1$ of the elliptic vertex algebras
$\elp$ are defined for any~$N$, extending the previously known case
of $N=2$. They realise deformed, centrally
extended double Yangian structures $\dey$. 
As in the quantum affine algebras $\uq{q}{N}$, and quantum elliptic
affine algebras $\elp$, these algebras contain
subalgebras at critical values of the central charge $c=-N-Mr$ 
($M$ integer, $2r=\ln p/\ln q$), which 
become Abelian when $c=-N$ or $2r=Nh$ for $h$ integer. Poisson
structures and quantum exchange relations are derived for their 
abstract generators.
\end{abstract}

\vfill
\vfill

\rightline{LAPTH-716/99}
\rightline{PAR-LPTHE 99-05}
\rightline{DTP-99-7}
\rightline{math.QA/9905100}
\rightline{February 1999}
\rightline{Revised: May 1999}

\newpage
\pagestyle{plain}
\setcounter{page}{1}

\section{Introduction}

Elliptic deformation of a finite Lie algebra was first proposed by
Sklyanin \cite{Skl}. The incorporation of a central extension leading
to the notion of vertex elliptic affine algebra $\elpa{2}$ was
achieved in \cite{FIJKMY} in the $R$-matrix formalism using
Belavin--Baxter 8-vertex matrix \cite{Ba}. Its quasi-Hopf structure as
a twist of $\uq{q}{2}$ was elucidated in \cite{Fronsdal}; a universal
formula for the twist was proposed in \cite{JKOS} and allowed to
define an $R$-matrix formulation for the general case $\elp$.

\medskip

Vertex elliptic affine algebras were shown to possess at least two
inequivalent limits with one parameter less. 
The limit $p \rightarrow 0$, together with a
suitable renormalisation of the generators, leads 
back to the $R$-matrix structure of the quantum group 
${\cal U}_{q}(\widehat{sl}(N)_{c})$ as shown in  
\cite{FIJKMY} for $N=2$.  

\medskip

Another limit was established for $N=2$ in \cite{KLP}, defining a so-called 
scaled limit ${\cal A}_{\hbar,\eta}(\widehat{sl}(2)_{c})$.  It was also 
described in \cite{JKM,Konno}.  It is obtained by sending $p,q$ and the 
spectral parameter $z$ to 1 while retaining as finite parameters $\ln p/\ln q$ 
and $\ln z/\ln q$.  A subsequent limit $\ln p/\ln q \rightarrow \infty$ leads 
\cite{KLP} to an algebraic structure of same $R$-matrix formulation as the 
centrally extended double Yangian ${\cal D}Y (sl(2))_{c}$ introduced in 
\cite{KT,Kh}.  Other definitions using differently
normalised $R$-matrices have been proposed for 
this object \cite{KoIo,Io}.  
The scaled algebra ${\cal A}_{\hbar,\eta} (\widehat{sl}(2)_{c})$ was
recently shown \cite{KLCP} to be realised by the screening currents
for the Sine--Gordon model at the free fermion point (see also
\cite{JKM,Konno}). At this time, only for $N=2$ has the scaling limit
been defined. 

\medskip 

It has been known for a while that quantum affine algebras
$\cU_q(\widehat{\cG})_c$ admit infinite dimensional centres at critical
values of 
the central charge \cite{RSTS}. Associated non-linear Poisson
structures realising $q$-deformed classical $\cW_N$ algebras were
derived for $\cG=sl(N)$ \cite{FR,AKOS} and quantised \cite{AKOS,FF}. A
similar, 
albeit more complex, structure of sub-exchange algebras, extended
centres or Abelian algebras, and classical (Poisson bracket) limits
thereof, was uncovered for $\elpa{2}$ \cite{AFRS1,AFRS2} and $\elp$
\cite{AFRS3,AFRS5}. They realise new quantum versions of the already
known $q$-$\cW_N$ Poisson structure \cite{FR,AKOS}. No such results are
known for the scaled algebras ${\cal A}_{\hbar,\eta}
(\widehat{sl}(2)_{c})$\footnote{This was pointed out to us by
  V.~Korepin and P.~Sorba.} and their yet to be defined generalisations
to $sl(N)$.

Our purpose therefore is to apply to the scaled limits 
${\cal  A}_{\hbar,\eta} (\widehat{sl}(N)_{c})$, once they are
consistently 
defined, the same systematic study undertaken in \cite{FR} for
affine quantum algebras and \cite{AFRS1,AFRS2,AFRS3,AFRS5} for vertex
elliptic affine algebras, i.e. deduce from the $R$-matrix formulation
the existence of extended centres (or at least Abelian subalgebras),
and compute the related Poisson structures and the related
quantisations. 

The ultimate purpose is to try to define new types of classical and
quantum deformed $\cW_N$-algebras. However a particular difficulty
arises here owing to the current state of understanding of scaled
algebras as univocally defined algebraic structures. 
Indeed, definition of the full 
algebraic structure requires a proper definition of the mode expansion
both for generating functionals and structure functions, and there may
exist several inequivalent expansions  
(and therefore inequivalent algebras) for one single abstract
structure.  
The structure functions here turn out to be single-period meromorphic
functions  
with an essential singularity at infinity.  It was commented in
\cite{KLP} that  
(in the rational limit) two different structures could be defined from one 
single $R$-matrix formulation, depending whether the modes of the
generators were defined to be 
continuous-labelled Fourier transforms with respect to the spectral parameter 
$\beta$, leading to the algebra ${\cal A}_{\hbar,0} (\widehat{sl}(2)_{c})$ 
\cite{KLP}, or discrete-labelled coefficient of power expansion in $\beta$, 
leading to the double Yangian ${\cal D}Y (sl(2))_{c}$ \cite{KT,Kh}.  
Similarly the deformed version of the exchange algebra, and the subsequent 
Poisson structures on the critical surfaces, will give rise to inequivalent 
algebraic structures depending upon the type of mode
expansion eventually  
chosen.  This discussion goes 
beyond the scope of the current presentation and will be left for
further studies.   
Identification of, and discussions on algebraic structures are  
thus generically understood here \emph{at the level of the generating
  functionals and $R$ matrix structures}.

\medskip

Within this framework we address the two basic problems of this paper as 
follows: 

In a first part, we define an $R$-matrix formulation for a scaling 
limit $q \rightarrow 1$ of the vertex algebra $\elp$.  We prove explicitly 
that the scaling limit structure matrix obeys the Yang--Baxter equation, 
unitarity, crossing-symmetry and quasi-periodicity relations.  We define a 
subsequent limit $\ln p/\ln q \rightarrow \infty$.  Its $R$-matrix
structure is  
strongly reminiscent of the double Yangian $\deysr$ in the 
formulation of \cite{KoIo,Io}, albeit with a different normalisation 
factor.  As a consequence the structure at $2r=\ln p/\ln q < \infty$ can be 
interpreted as a deformation of the double Yangian of which more will
be said in the conclusion, hence it will be denoted 
$\dey$.  

In a second part, we show that the scaled algebra $\dey$ contains 
non-linear exchange subalgebras when $c=-N-Mr$ for $M \in \ZZ$.  When $M=0$ 
(that is $c=-N$) these algebras are Abelian.  When $M \ne 0$, the 
algebras become Abelian when $2r = Nh = 2 r_{crit}$ for $h \in \ZZ
\backslash \{0\}$.   
Non-linear Poisson structures are then derived on these Abelian 
subalgebras, and the exchange algebras at $M \ne 0$ can be interpreted as 
consistent quantisations inside $\dey$ of the Poisson structures with 
$\hbar \sim 2(r-r_{crit})$.  They contain higher-spin operators obtained as 
before \cite{AFRS5} 
by the $R$-matrix ``twisted'' trace formula of the Lax operators.  
It appears here that the existence of exchange subalgebras on critical 
surfaces and Abelian limits thereof on critical subalgebras, survive 
the scaling limit simply by taking directly the same scaling limit in
the characteristic relations and the structure functions.
Interpretations of the Poisson algebras and exchange  
algebras, in connection with Virasoro and ${\cal W}_{N}$ algebras, require 
to establish the exact mode expansion to be used on these generating 
functionals, and will be left for further studies.

In conclusion we indicate several directions of interest which we mean
to pursue in the next future. 

\section{Degeneration limit of the elliptic algebra $\elp$}
\setcounter{equation}{0}

\subsection{Definition of the $sl(N)$ elliptic $R$-matrix}

The $sl(N)$ elliptic $R$-matrix in $\mathrm{End}(\CC^N) \otimes 
\mathrm{End}(\CC^N)$, associated to the $\ZZ_{N}$-vertex model, is defined 
as follows \cite{Bela,ChCh}:
\begin{equation}
  R(z,q,p) = z^{2/N-2} \frac{1}{\kappa(z^2)} 
  \frac{\vartheta\car{\half}{\half}(\zeta,\tau)} 
  {\vartheta\car{\half}{\half}(\xi+\zeta,\tau)} \,\, 
  \sum_{(\alpha_1,\alpha_2)\in\ZZ_N\times\ZZ_N} 
  W_{(\alpha_1,\alpha_2)}(\xi,\zeta,\tau) \,\, I_{(\alpha_1,\alpha_2)} 
  \otimes I_{(\alpha_1,\alpha_2)}^{-1} \,,
  \label{eq27}
\end{equation}
where the variables $z,q,p$ are related to the variables 
$\xi,\zeta,\tau$ by
\begin{equation}
  z=e^{i\pi\xi} \,,\qquad q=e^{i\pi\zeta} \,,\qquad p=e^{2i\pi\tau} \,.
  \label{eq28}
\end{equation}
The Jacobi theta functions with rational characteristics $\gamma = 
(\gamma_1,\gamma_2) \in \sfrac{1}{N} \ZZ \times \sfrac{1}{N} \ZZ$ are 
defined by:
\begin{equation}
  \vartheta\car{\gamma_1}{\gamma_2}(\xi,\tau) = \sum_{m \in \ZZ}
  \exp\Big(i\pi(m+\gamma_1)^2\tau + 2i\pi(m+\gamma_1)(\xi+\gamma_2)
  \Big) \,.
  \label{eqvartheta}
\end{equation}
The normalisation factor is chosen as follows:
\begin{equation}
  \frac{1}{\kappa(z^2)} = \frac{(q^{2N}z^{-2};p,q^{2N})_\infty
    \, (q^2z^2;p,q^{2N})_\infty \, (pz^{-2};p,q^{2N})_\infty \,
    (pq^{2N-2}z^2;p,q^{2N})_\infty} {(q^{2N}z^2;p,q^{2N})_\infty
    \, (q^2z^{-2};p,q^{2N})_\infty \, (pz^2;p,q^{2N})_\infty \,
    (pq^{2N-2}z^{-2};p,q^{2N})_\infty} \,.
  \label{eq29}
\end{equation}
The functions $W_{(\alpha_1,\alpha_2)}$ are given by
\begin{equation}
  W_{(\alpha_1,\alpha_2)}(\xi,\zeta,\tau) = 
  \frac{\vartheta\car{\half+\alpha_1/N} {\half+\alpha_2/N}(\xi+\zeta/N,\tau)} 
  {N\vartheta\car{\half+\alpha_1/N} {\half+\alpha_2/N}(\zeta/N,\tau)} \,.
  \label{eq211}
\end{equation}
The matrices $I_{(\alpha_1,\alpha_2)}$ are defined by
\begin{equation}
  I_{(\alpha_1,\alpha_2)} = g^{\alpha_2} \, h^{\alpha_1} \,,
  \label{eq212}
\end{equation}
where the $N \times N$ matrices $g$ and $h$ are given by $g_{ij} =
\omega^i\delta_{ij}$ and $h_{ij} = \delta_{i+1,j}$, the addition of 
indices being understood modulo $N$ and $\omega = e^{2i\pi/N}$.  \\
Let us set
\begin{equation}
  S(\xi,\zeta,\tau) = \sum_{(\alpha_1,\alpha_2)\in\ZZ_N\times\ZZ_N}
  W_{(\alpha_1,\alpha_2)}(\xi,\zeta,\tau) \,\, I_{(\alpha_1,\alpha_2)}
  \otimes I_{(\alpha_1,\alpha_2)}^{-1} \,.
  \label{eq210}
\end{equation}
The matrix $S$ is $\ZZ_N$-symmetric, that is
\begin{equation}
  S_{a+s\,,\,b+s}^{c+s\,,\,d+s} = S_{a\,,\,b}^{c\,,\,d}
  \label{eq213}
\end{equation}
for any indices $a,b,c,d,s \in \ZZ_N$ (the addition of indices being
understood modulo $N$) and the non-vanishing elements of the matrix
$S$ are of the type $S_{a\,,\,b}^{c\,,\,a+b-c}$.  One finds explicitly:
\begin{equation}
  S_{a\,,\,b}^{c\,,\,a+b-c}(\xi,\zeta,\tau) = 
  \frac{\vartheta\car{(b-a)/N+\half}{\half}(\xi+\zeta,N\tau)} 
  {\vartheta\car{(c-a)/N+\half}{\half}(\zeta,N\tau)} \,\, 
  \frac{\displaystyle\prod_{k=0,k \ne b-c}^{N-1} 
    \vartheta\car{k/N+\half}{\half}(\xi,N\tau)} 
  {\displaystyle\prod_{k=1}^{N-1}\vartheta\car{k/N+\half} {\half}(0,N\tau)} \,.
  \label{eq215}
\end{equation}
Using the relation
\begin{equation}
  \prod_{k=1}^{N-1}\vartheta\car{k/N+\half} {\half}(\xi,N\tau) 
  = p^{(N-1)(N-2)/24} \; \frac{\displaystyle\left[ \prod_{m=1}^\infty 
    (1-p^{Nm}) \right]^{N}}{\displaystyle\prod_{m=1}^\infty (1-p^{m})} \; 
  \frac{\Theta_{p}(z^2)}{\Theta_{p^N}(z^2)} \;,
\end{equation}
one finds the following expression in terms of the Jacobi $\Theta$ functions 
for the elements $R_{a\,,\,b}^{c\,,\,a+b-c}$ of the matrix $R(z,q,p)$:
\begin{eqnarray}
  R_{a\,,\,b}^{c\,,\,a+b-c} &=& \frac{1}{\kappa(z^2)} \,\,
  z^{\frac{2}{N}\,(1-a)}  
  \, q^{\frac{2b}{N}} \, p^{-\frac{(a-c)(b-c)}{N}} \;
  \frac{\displaystyle\left[  
    \prod_{m=1}^\infty (1-p^{Nm}) \right]^{3}}{\displaystyle\left[ 
    \prod_{m=1}^\infty (1-p^{m}) \right]^{3}} \nonumber \\
  && \qquad \qquad \times \frac{\Theta_{p^N}(p^{N+b-a}q^2z^2)} 
  {\Theta_{p^N}(p^{N+b-c}z^2)\,\Theta_{p^N}(p^{N-a+c}q^2)} 
  \frac{\Theta_p(q^2)\,\Theta_p(pz^2)}{\Theta_p(q^2z^2)} \,.
  \label{eq219}
\end{eqnarray}

\subsection{Scaling limit of the $sl(N)$ elliptic $R$-matrix}

The scaling limit of the $R$ matrix (\ref{eq219}) is defined by
\begin{equation}
  q \rightarrow 1 \,, \qquad z = q^{i\beta/\pi} \,, \qquad p = q^{2r} \,, 
  \qquad \mbox{where $\beta,r$ are kept fixed} \,.
  \label{eqscal}
\end{equation}
Evaluation of the limit of the elements (\ref{eq219}) requires regularisation 
of potentially divergent terms (for instance by exchanging limits and infinite 
products symbols).  This in turn implies that the crucial identities (such as 
Yang--Baxter equation, unitarity, \ldots) of the scaling limit of the $R$ 
matrix be checked again (see Theorem \ref{thmone}).  \\
One finds in the limit (\ref{eqscal}) that
\begin{equation}
  \frac{\Theta_{p^n}(p^u q^{2v} z^{2m})}
  {\Theta_{p^{n'}}(p^{u'} q^{2v'} z^{2m'})} \quad \rightarrow \quad 
  \frac{\sin\frac{\pi}{nr}(ur+v+im\beta/\pi)}
  {\sin\frac{\pi}{n'r}(u'r+v'+im'\beta/\pi)} \;.
\end{equation}
It follows that in the scaling limit the matrix elements (\ref{eq219}) 
become
\begin{equation}
  R_{a\,,\,b}^{c\,,\,a+b-c}(\beta,r) = -\frac{1}{N} \; S_{0}(\beta) \; 
  \frac{\displaystyle \sin\frac{i\beta+\pi+(b-a)\pi r}{Nr}} 
  {\displaystyle \sin\frac{i\beta+(b-c)\pi r}{Nr}
    \,\sin\frac{\pi-(a-c)\pi r}{Nr}} \;  
  \frac{\displaystyle \sin\frac{i\beta}{r}\,\sin\frac{\pi}{r}}
  {\displaystyle \sin\frac{i\beta+\pi}{r}} \;.
  \label{eq219scal}
\end{equation}
Remark that the factor $1/N$ arises from a regularisation of the limit of the 
ratio of the infinite product over $m$ in (\ref{eq219}). More
precisely, the infinite product arises from the $\xi \rightarrow 0$
limit of a ratio of theta functions. The regularisation requires to
exchange the limits $\xi \rightarrow 0$ with the limit $q\rightarrow 1$
in this ratio of theta functions. \\
The normalisation factor $S_{0}(\beta)$ is defined by
\begin{equation}
  S_{0}(\beta) = \displaystyle 
  \frac{S_{2}(-\frac{i\beta}{\pi} \vert r,N) \, 
    S_{2}(1+\frac{i\beta}{\pi} \vert r,N)}
  {S_{2}(\frac{i\beta}{\pi} \vert r,N) \, 
    S_{2}(1-\frac{i\beta}{\pi} \vert r,N)} \,,
  \label{eqso}
\end{equation}
where $S_{2}(x \vert \omega_{1},\omega_{2})$ is the Barnes' double sine 
function \cite{JiMi} of periods $\omega_{1}$ and $\omega_{2}$ defined by
\begin{equation}
  S_{2}(x \vert \omega_{1},\omega_{2}) = 
  \frac{\Gamma_{2}(\omega_{1}+\omega_{2}-x \vert \omega_{1},\omega_{2})}
  {\Gamma_{2}(x \vert \omega_{1},\omega_{2})} \;,
\end{equation}
where
\begin{equation}
  \Gamma_{2}(x \vert \omega_{1},\omega_{2}) = \exp\left( 
  \frac{\partial}{\partial s} \sum_{n_{1},n_{2} \ge 0} 
  (x+n_{1}\omega_{1}+n_{2}\omega_{2})^{-s}\Bigg\vert_{s=0} \right) \,.
\end{equation}
It satisfies
\begin{equation}
  \frac{S_{2}(x+\omega_{1} \vert \omega_{1},\omega_{2})}{S_{2}(x \vert 
    \omega_{1},\omega_{2})} = \frac{1}{2\sin\frac{\pi x}{\omega_{2}}} \,,
  \qquad \qquad
  \frac{S_{2}(x+\omega_{2} \vert \omega_{1},\omega_{2})}{S_{2}(x \vert 
    \omega_{1},\omega_{2})} = \frac{1}{2\sin\frac{\pi x}{\omega_{1}}} \,.
\end{equation}
$-S_{0}(\beta)$ is the limit of the normalisation factor $1/\kappa(z^2)$ of 
the matrix (\ref{eq27}).

\medskip

\noindent In the particular case $N=2$, one recovers explicitly \cite{Konno}:
\begin{equation}
  R(\beta,r) = - S_{0}(\beta) \left( 
  \begin{array}{cccc} 
    \displaystyle \frac{\cosh\frac{\beta}{2r} \; \cosh\frac{i\pi}{2r}} 
    {\cosh\frac{i\pi-\beta}{2r}} & 0 & 0 & \displaystyle 
    -\frac{\sinh\frac{\beta}{2r} \; \sinh\frac{i\pi}{2r}} 
    {\cosh\frac{i\pi-\beta}{2r}} \\
    0 & \displaystyle -\frac{\sinh\frac{\beta}{2r} \; 
      \cosh\frac{i\pi}{2r}} {\sinh\frac{i\pi-\beta}{2r}} & \displaystyle 
    \frac{\cosh\frac{\beta}{2r} \; \sinh\frac{i\pi}{2r}} 
    {\sinh\frac{i\pi-\beta}{2r}} & 0 \\
    0 & \displaystyle \frac{\cosh\frac{\beta}{2r} \; 
      \sinh\frac{i\pi}{2r}} {\sinh\frac{i\pi-\beta}{2r}} & \displaystyle 
    -\frac{\sinh\frac{\beta}{2r} \; \cosh\frac{i\pi}{2r}} 
    {\sinh\frac{i\pi-\beta}{2r}} & 0 \\
    \displaystyle -\frac{\sinh\frac{\beta}{2r} \; \sinh\frac{i\pi}{2r}} 
    {\cosh\frac{i\pi-\beta}{2r}} & 0 & 0 & \displaystyle 
    \frac{\cosh\frac{\beta}{2r} \; \cosh\frac{i\pi}{2r}} 
    {\cosh\frac{i\pi-\beta}{2r}} \\
  \end{array} \right) \,.  \label{eqmatdy}
\end{equation}

\bigskip

\begin{thm}\label{thmone}
  The matrix $R(\beta,r)$ satisfies the following properties: \\
  -- Yang--Baxter equation:
  \begin{equation}
    R_{12}(\beta_{1}-\beta_{2}) \, R_{13}(\beta_{1}-\beta_{3}) \, 
    R_{23}(\beta_{2}-\beta_{3}) = R_{23}(\beta_{2}-\beta_{3}) \, 
    R_{13}(\beta_{1}-\beta_{3}) \, R_{12}(\beta_{1}-\beta_{2}) \,,
    \label{eq221}
  \end{equation}
  -- Unitarity:
  \begin{equation}
    R_{12}(\beta) \, R_{21}(-\beta) = 1 \,,
    \label{eq222}
  \end{equation}
  -- Crossing-symmetry:
  \begin{equation}
    R_{12}(\beta)^{t_2} \, R_{21}(-\beta+Ni\pi)^{t_2} = 1 \,,
    \label{eq223}
  \end{equation}
  -- Quasi-periodicity:
  \begin{equation}
    {\widehat R}_{12}(\beta-i\pi r) = (h \otimes \II)^{-1} \, {\widehat 
      R}_{21}(-\beta)^{-1} \, (h \otimes \II) \,,
    \label{eq225}
  \end{equation}
  where
  \begin{equation}
    {\widehat R}_{12}(\beta) = 
    \frac{\displaystyle \sin\frac{\pi-i\beta}{N}}
    {\displaystyle \sin\frac{i\beta}{N}} \, R_{12}(\beta) \,.
    \label{eq226}
  \end{equation}
\end{thm}
Note that the antisymmetry property of the elliptic $R$ matrix (\ref{eq27}) 
under $z \rightarrow -z$ does not make sense in the scaling limit 
(\ref{eqscal}) for obvious reasons.

\medskip

\noindent \textbf{Proof:} The proof of Theorem \ref{thmone} is based on 
explicit computations using the form of the matrix (\ref{eq219scal}) and 
standard arguments about meromorphic functions.  Illustration is given on the 
example of the unitarity relation.  We recall that the explicit form for the 
$R$ matrix is:
\begin{equation}
  R_{a\,,\,b}^{c\,,\,a+b-c}(\beta,r) = -\frac{1}{N} \; 
  S_{0}(\beta) \; \displaystyle  
  \frac{\sin\frac{i\beta}{r}\,\sin\frac{\pi}{r}} 
  {\sin\frac{i\beta+\pi}{r}} \;
  \left( \cotan\frac{i\beta+(b-c)\pi r}{Nr} -
    \cotan\frac{-\pi+(a-c)\pi r}{Nr} 
  \right) \;.
\end{equation}
Denoting $x = -\pi/Nr$ and $y = i\beta/Nr$, eq.  (\ref{eq222}) reads:
\begin{eqnarray}
  R_{12}(x,y) R_{21}(x,-y) &=& {\cal N} \sum_{a,b,c=1}^N \sum_{i=1}^N 
  \left(
    \Big(\cotan(y+\frac{(b-i)\pi}{N}) - \cotan(x+\frac{(a-i)\pi}{N})\Big) 
  \right. \nonumber \\
  &&
  \left.
    \Big(-\cotan(x+\frac{(c-i)\pi}{N}) + \cotan(-y+\frac{(i-a-b+c)\pi}{N})\Big)
  \right) 
  e_{a}^c \otimes e_{b}^{a+b-c} \nonumber \;, \\
  \label{equnit}
\end{eqnarray}
where ${\cal N}^{-1} = N^2(\cotan^2Nx - \cotan^2Ny)$; this overall factor 
comes from the product of the index-independent factors $-\frac{1}{N} \; 
S_{0}(\beta) \; \displaystyle 
\frac{\sin\frac{i\beta}{r}\,\sin\frac{\pi}{r}} 
{\sin\frac{i\beta+\pi}{r}}$.  Notice in particular that 
$S_{0}(\beta)S_{0}(-\beta) = 1$ as obviously follows from (\ref{eqso}).

\medskip

\noindent We now prove that the coefficient function in (\ref{equnit}) is 
proportional to $\delta_{ac}$ by analyzing its poles as a function 
of $y$.

\medskip

\noindent If $a \ne c$, this function has simple poles at $y = k\pi/N$ for all 
integer $k$ since $a$, $b$, $c$ and $i$ are integers in (\ref{equnit}).
One immediately sees on its explicit expression that the residue  
is in fact zero.  Hence, the function does not depend on $y$.  One then 
evaluates it at $y = x + (a-b)\pi/N$ where it becomes zero. 

\medskip

\noindent If now $a=c$, this function can be directly computed.  The relation 
(\ref{equnit}) becomes:
\begin{equation}
  R_{12}(x,y) R_{21}(x,-y) = {\cal N} \sum_{a,b=1}^N \sum_{i=1}^N 
  \left(
    \cotan^2(x+\frac{(a-i)\pi}{N}) - \cotan^2(y+\frac{(b-i)\pi}{N}) 
  \right)
  e_{a}^a \otimes e_{b}^b \;.
\label{eqcotan}
\end{equation}
Since the sum over $i$ is manifestly cyclic, the expression does not depend on 
$a$ and $b$.  It is then explicitly evaluated from the equality \cite{Prud}
\begin{equation}
  \sum_{i} \cotan^2(x+i\pi/N) = N^2 - N + N^2 \cotan^2Nx \;.
  \label{eqprud}
\end{equation}
Combining (\ref{eqcotan}), (\ref{eqprud}) and the expression of the 
normalisation factor ${\cal N}$, one gets the desired result.

\bigskip

\noindent The crossing-symmetry (\ref{eq223}) is proved by similar arguments.  
The quasi-periodicity relation (\ref{eq225}) is essentially trivial to prove.  
The Yang--Baxter equation by contrast requires a careful analysis of the 
behaviour at the poles, involving as it does combinations of double and simple 
poles.
\finproof

\bigskip

\textbf{Remark:} Using the crossing symmetry and the unitarity
properties of $R_{12}$, one can exchange inversion and transposition
for the matrix $R_{12}$ as follows (the same property also holds for
the matrix $\widehat R_{12}$):
\begin{equation}
  \left( R_{12}(\beta)^{t_2}
  \right)^{-1}
  =
  \left( R_{12}(\beta - Ni\pi)^{-1}
  \right)^{t_2} \;.
  \label{eq:crossunit}
\end{equation}

\subsection{Definition of the algebra}

The limit of the elliptic quantum algebra structure ${\cal 
A}_{q,p}(\widehat{sl}(N)_{c})$ defined by the scaling procedure (\ref{eqscal}) 
leads to an algebraic structure parametrised by the limit 
$R$-matrix for abstract generators $L(\beta)$ encapsulated into a $N \times N$ 
matrix
\begin{equation}
  L(\beta) = \left(\begin{array}{ccc} L_{11}(\beta) & \cdots &
  L_{1N}(\beta) \cr \vdots && \vdots \cr L_{N1}(\beta) & \cdots &
  L_{NN}(\beta) \cr \end{array}\right) \,.
\end{equation}
One defines $\deygl$ by imposing the following 
constraints on the $L_{ij}(\beta)$ (with the matrix $\widehat R_{12}$ given by 
eq.  (\ref{eq226})):
\begin{equation}
  \widehat R_{12}(\beta_{1}-\beta_{2}) \, L_1(\beta_{1}) \, L_2(\beta_{2}) = 
  L_2(\beta_{2}) \, L_1(\beta_{1}) \, \widehat 
  R_{12}^{*}(\beta_{1}-\beta_{2}) \,,
  \label{eq232}
\end{equation}
where $L_1(\beta) \equiv L(\beta) \otimes \II$, $L_2(\beta) \equiv \II \otimes 
L(\beta)$ and $\widehat R^{*}_{12}$ is defined by $\widehat 
R^{*}_{12}(\beta,r) \equiv \widehat R_{12}(\beta,r-c)$.  
\\
The matrix 
$\widehat R^{*}_{12}$ obeys also the unitarity, crossing-symmetry and 
quasi-periodicity conditions of Theorem~\ref{thmone}.  \\
The $q$-determinant $q$-$\det L(\beta)$ is given by ($\eps(\sigma)$ being the 
signature of the permutation $\sigma$)
\begin{equation}
  q\mbox{-}\det L(\beta) \equiv \sum_{\sigma\in{\mathfrak S}_N}
  \eps(\sigma) \prod_{k=1}^N L_{k,\sigma(k)}(\beta -i\pi(k-N-1)) \;.
  \label{eq233}
\end{equation}
It lies in the centre of $\deygl$; it can be
``factored  out'' and set to the value $1$ so as to get
\begin{equation}
  \dey = \deygl/ 
  \langle q\mbox{-}\det L - 1 \rangle \,.
  \label{eq234}
\end{equation}
By analogy with the full elliptic case we now introduce the following
two matrices: 
\begin{subeqnarray}
  \label{eq235}
  && L^+(\beta) \equiv L(\beta-\half\,i\pi c) \,, \\
  && L^-(\beta) \equiv h \, L(\beta-i\pi r) \, h^{-1} 
  \,.
\end{subeqnarray}
They obey coupled exchange relations following from (\ref{eq232}) and
periodicity/unitarity properties of the matrices
$\widehat R_{12}$ and $\widehat R^{*}_{12}$:
\begin{subeqnarray}
  && \widehat R_{12}(\beta_{1}-\beta_{2}) \, L^\pm_1(\beta_{1}) \,
  L^\pm_2(\beta_{2}) = L^\pm_2(\beta_{2}) \, L^\pm_1(\beta_{1}) \,
  \widehat R^{*}_{12}(\beta_{1}-\beta_{2}) \,, \\
  && \widehat R_{12}(\beta_{1}-\beta_{2}-\half\,i\pi c) \,
  L^+_1(\beta_{1}) \, L^-_2(\beta_{2}) = L^-_2(\beta_{2}) \, L^+_1(\beta_{1})
  \, \widehat R^{*}_{12}(\beta_{1}-\beta_{2}+\half\,i\pi c) \,.
  \label{eq236}
\end{subeqnarray}

\paragraph{Important remark:} The algebraic structure (\ref{eq236})
and the properties (\ref{eq221})--(\ref{eq225}) and
(\ref{eq:crossunit}) are identical to the algebraic structure and
properties of $R$-matrices for the full elliptic case, provided that
the parameters $p$, $z$ be replaced by their scaling limits $r$,
$\beta$; the multiplicative shifts are replaced by additive shifts;
and the scalar factor in $\widehat R$ be replaced by its
renormalised scaling limit (\ref{eq226}). It follows that all
conclusions touching the existence and structure functions of
algebraic substructures (exchange algebras; Abelian algebras and
Poisson algebras) can be transposed from the elliptic case to the
scaling limit by suitable changes of the structure functions,
provided that their proofs in \cite{AFRS3,AFRS5} use only the
relations equivalent to (\ref{eq221}--\ref{eq225}),
(\ref{eq:crossunit}) and (\ref{eq236}).

\section{Exchange algebras on the double Yangian $\dey$}
\setcounter{equation}{0}

\begin{lemm}
  Let 
  $\cL^{(M)}(\beta) = ( h^{-M} L^+(\beta-\half\,i\pi c))^t
  \; \wL^-(\beta)$, 
  where $\wL^\pm(\beta) \equiv \left(L^\pm(\beta)^t \right)^{-1}$ and
  $M\in\ZZ$.  
  The operators $\cL(\beta)$ have the following exchange properties with
  the generators $L$ of $\dey$ when the parameters $c$ and $r$ obey 
  the relation $c= -N - Mr$:
  \begin{equation}
    \cL^{(M)}_1(\beta_1) L_2(\beta_2) = 
    F(M,\beta_2-\beta_1)
    L_2(\beta_2) 
    \left(\widehat R_{21}^*(\beta_2-\beta_1 + i\pi c)^{-1}
    \right)^{t_1}
    \cL^{(M)}_1(\beta_1)
    \widehat R_{21}^*(\beta_2 - \beta_1 + i\pi c)^{t_1}   \;,
    \label{eq31}
  \end{equation}
  where 
  \begin{equation}
    F(M,\beta) = 
    \left\{ 
      \begin{array}{ll}
        \displaystyle \prod_{k=0}^{M-1}
        \frac{\displaystyle \sin^2\frac{i\beta+k\pi r}{N}}
        {\displaystyle 
          \sin\frac{i\beta+k\pi r+\pi}{N} \; 
          \sin\frac{i\beta+k\pi r-\pi}{N}} 
        & \mbox{for $M>0$} \,, \\ \\
        \displaystyle \prod_{k=1}^{|M|}
        \frac{\displaystyle \sin\frac{-i\beta+k\pi r+\pi}{N} \;
          \sin\frac{-i\beta+k\pi r-\pi}{N}} 
        {\displaystyle \sin^2\frac{-i\beta+k\pi r}{N}}
        & \mbox{for $M<0$} \,, \\
      \end{array} 
    \right.
    \label{eq:FM}
  \end{equation}
  and 
  \begin{equation}
    F(0,\beta) = 1 \;.
  \end{equation}
  \label{lem1}
\end{lemm}

\noindent
We give here an extensive proof of Lemma \ref{lem1} as an illustration
of the transposition procedure given in the previous remark. Other
proof will simply call back to this key feature. 

\noindent
\textbf{Proof:}  
Let us prove the lemma in terms of $L^+(\beta)$.
One has
\begin{equation}
  \cL^{(M)}_1(\beta_1) L_2^+(\beta_2) = L_1^+(\beta_1-\half\,i\pi
  c)^{t_1} (h^{-M}_1)^{t_1} \wL_1^-(\beta_1) L^+_2(\beta_2) \;.
\end{equation}
Exchanging $\cL^{(M)}_1(\beta_1)$ with $L^+_2(\beta_2)$ requires the
following  
relations obtained from (\ref{eq236}):
\begin{subeqnarray}
  \left(\widehat R_{21}(\beta_2-\beta_1 - \half \,  i\pi c)^{t_1}
  \right)^{-1} \, L_2^+(\beta_2) \, \widetilde L_1^-(\beta_1) 
  &=&
  \widetilde L_1^-(\beta_1) \, L_2^+(\beta_2) \,
  \left(\widehat R_{21}^*(\beta_2-\beta_1 + \half \, i\pi c)^{t_1}
  \right)^{-1} 
  \slabel{eq32a}
  \\
  L_1^+(\beta_1)^{t_1} 
  \left( \widehat R_{21}(\beta_2-\beta_1)^{-1}
  \right)^{t_1} L_2^+(\beta_2) 
  &=&
  L_2^+(\beta_2) 
  \left( \widehat R_{21}^* (\beta_2-\beta_1)^{-1}
  \right)^{t_1} 
  L_1^+(\beta_1)^{t_1} \;.
  \slabel{eq32b}
\end{subeqnarray}
Using (\ref{eq32a}), one gets
\begin{eqnarray}
  \cL^{(M)}_1(\beta_1) L_2^+(\beta_2) &=&   
  L_1^+(\beta_1-\half\,i\pi
  c)^{t_1} 
  (h^{-M}_1)^{t_1} 
  \left(\widehat R_{21}(\beta_2-\beta_1 - \half \,  i\pi c)^{t_1}
  \right)^{-1}
  \nonumber\\ && \qquad\qquad
  L_2^+(\beta_2) 
  \wL_1^-(\beta_1)
  \widehat R_{21}^*(\beta_2-\beta_1 + \half \, i\pi c)^{t_1} 
  \nonumber\\
  &=&
  L_1^+(\beta_1-\half\,i\pi
  c)^{t_1} 
  (h^{-M}_1)^{t_1} 
  \left(\widehat R_{21}(\beta_2-\beta_1 - \half \,  i\pi c)^{t_1}
  \right)^{-1}
  \nonumber\\ && \qquad \qquad
  (h^{M}_1)^{t_1}  L_2^+(\beta_2)   (h^{-M}_1)^{t_1} 
  \wL_1^-(\beta_1)
  \widehat R_{21}^*(\beta_2-\beta_1 + \half \, i\pi c)^{t_1} \;.
  \label{eq33}
\end{eqnarray}
{}From the crossing unitarity property (\ref{eq:crossunit}) of
$\widehat R_{21}$, this equation can be written as
\begin{eqnarray}
  \cL^{(M)}_1(\beta_1) L_2^+(\beta_2) &=&   L_1^+(\beta_1-\half\,i\pi
  c)^{t_1} 
  \left(h^M_1 \widehat R_{21}(\beta_2-\beta_1 - \half \,  i\pi c
    -Ni\pi)^{-1} h^{-M}_1  
  \right)^{t_1}  
  \nonumber \\ && \qquad \qquad
  L_2^+(\beta_2) 
  (h^{-M}_1)^{t_1}
  \wL_1^-(\beta_1)
  \widehat R_{21}^*(\beta_2-\beta_1 + \half \, i\pi c)^{t_1} \;.
  \label{eq34}
\end{eqnarray}
Now sufficient conditions that allow substitution of eq. (\ref{eq32b})
into (\ref{eq34}) are 
\begin{equation}
  M r = -c - N \qquad \mbox{where} \qquad M\in\ZZ \,.
  \label{eq:surf}
\end{equation}
At this point, it is necessary to distinguish whether $M=0$ or not. 

\medskip

\begin{itemize}
\item 
In the case $M=0$, the condition (\ref{eq:surf}) fixes the central
charge $c$ to the critical value~$-N$. The argument of
$\widehat R_{21}$ in eq. (\ref{eq34}) is then exactly the one needed
to use eq. (\ref{eq32b}). Hence inserting (\ref{eq32b}) into
(\ref{eq34}), one obtains (with $\cL(\beta)\equiv  \cL^{(0)}(\beta)$)
\begin{eqnarray}
  && \cL_1(\beta_1)  L_2^+(\beta_2) = 
  \nonumber \\
  && \qquad = 
  L_2^+(\beta_2) 
  \left(\widehat R_{21}^* (\beta_2-\beta_1 + \half \, i\pi c )^{-1}
  \right)^{t_1}
  L_1^+(\beta_1-\half\,i\pi c)^{t_1} 
  \wL_1^-(\beta_1)
  \widehat R_{21}^*(\beta_2-\beta_1 + \half \, i\pi c)^{t_1} 
  \nonumber  \\
  && \qquad =
  L_2^+(\beta_2) 
  \left(\widehat R_{21}^* (\beta_2-\beta_1 + \half \, i\pi c )^{-1}
  \right)^{t_1}
  \cL_1(\beta_1) 
  \widehat R_{21}^*(\beta_2-\beta_1 + \half \, i\pi c)^{t_1} \;,
  \label{eq35}
\end{eqnarray}
which is formula (\ref{eq31}).

This proof reproduces exactly the proof given in \cite{AFRS3,AFRS5}
for the identification of exchange properties in the elliptic
case. The relation $c=-N-Mr$ is the scaling limit of the relation
$p^{M/2}=q^{-c-N}$ in \cite{AFRS5}. Identification of the relations
(\ref{eq221})-(\ref{eq225}) with scaling limits of \cite{AFRS5} is
crucial in this procedure. 

\item 
In the case $M \ne 0$, the argument of $\widehat R_{21}$ in
eq. (\ref{eq34}) is  $\beta_2-\beta_1+ \half \, i\pi c + i\pi Mr$. 
The use of the condition (\ref{eq:surf}) thus relies
on the quasiperiodicity property of $\widehat R_{21}$. 
More precisely, one has 
\begin{equation}
  h_1^M \widehat R_{21}(\beta + i\pi Mr) h_1^{-M}= F(-M,\beta) \widehat
  R_{21}(\beta) \;.
\end{equation}
Then on the line defined by (\ref{eq:surf}), the equation (\ref{eq34})
becomes 
\begin{eqnarray}
  \cL^{(M)}_1(\beta_1) L_2^+(\beta_2) &=&   
  F(-M,\beta_2-\beta_1+\half \, i\pi c)^{-1}  
  L_1^+(\beta_1-\half\,i\pi   c)^{t_1} 
  \left(\widehat R_{21}(\beta_2-\beta_1 + \half \,  i\pi c )^{-1}
  \right)^{t_1}
  \nonumber \\
  &&
  L_2^+(\beta_2) 
  (h^{-M}_1)^{t_1}
  \wL_1^-(\beta_1)
  \widehat R_{21}^*(\beta_2-\beta_1 + \half \, i\pi c)^{t_1} 
  \nonumber \\
  &=&
  F(-M,\beta_2-\beta_1+\half \, i\pi c)^{-1}  
  L_2^+(\beta_2) 
  \left(\widehat R_{21}^*(\beta_2-\beta_1 + \half \,  i\pi c )^{-1}
  \right)^{t_1}
  \nonumber \\
  &&
  L_1^+(\beta_1-\half\,i\pi   c)^{t_1} 
  (h^{-M}_1)^{t_1}
  \wL_1^-(\beta_1)
  \widehat R_{21}^*(\beta_2-\beta_1 + \half \, i\pi c)^{t_1} \;,
\end{eqnarray}
the last equality following from eq. (\ref{eq32b}). 
Hence we get, with $\beta=\beta_2-\beta_1 + \half \,  i\pi c $,
\begin{equation}
  \cL^{(M)}_1(\beta_1) L_2^+(\beta_2) = 
  F(-M,\beta)^{-1}  
  L_2^+(\beta_2) 
  \left(\widehat R_{21}^*(\beta)^{-1}
  \right)^{t_1}
  \cL^{(M)}_1(\beta_1)
  \widehat R_{21}^*(\beta)^{t_1}   \;.
  \label{eq36}
\end{equation}
Finally, one can check that $F(-M,\beta)^{-1} = F(M,\beta-i\pi c)$,
which achieves the proof.
\end{itemize}
\finproof

In the case $N=2$, the formula (\ref{eq:FM}) reduces to 
\begin{equation}
  F(M,\beta) = 
  \left\{ 
    \begin{array}{ll}
      (-1)^M \displaystyle \prod_{k=0}^{M-1}
      \tan^2\frac{i\beta+k\pi r}{2}
      & \mbox{for $M>0$} \,, \\ \\
      (-1)^M \displaystyle \prod_{k=1}^{|M|}
      \cotan^2\frac{i\beta-k\pi r}{2}
      & \mbox{for $M<0$} \,. \\
    \end{array} 
  \right.
  \label{eq:FM2}
\end{equation}

\subsection{Centre of $\dey$ at the critical level $c=-N$}
Using the previous lemma, it is straightforward to prove the following
theorem: 
\begin{thm}
  At the critical level $c=-N$, the operators generated by
  \begin{equation}
    t(\beta) \equiv {\rm Tr}[{\cal L}(\beta)] \equiv 
    {\rm Tr}\Big(L^+(\beta-\half\,i\pi c)^t \, \wL^-(\beta) \Big) 
    \label{eq:thm}
  \end{equation}
  lie in the centre of the algebra $\dey$.  
  \label{thm2}
\end{thm}

\medskip

\noindent\textbf{Proof:}
Equation (\ref{eq31}) leads to 
\begin{eqnarray}
  t(\beta_1) \, L_2(\beta_2)  
  &=&  
  {\rm Tr_1}[{\cL_1}(\beta_1)] \; L_2(\beta_2) \nonumber\\
  &=&
  L_2(\beta_2) \; {\rm Tr}_1
  \left[
  \left(\widehat R_{21}^* (\beta_2-\beta_1 + i\pi c )^{-1}
  \right)^{t_1}
  \cL_1(\beta_1) 
  \widehat R_{21}^*(\beta_2-\beta_1 + i\pi c)^{t_1} 
  \right] \;.
  \label{eq35bis}
\end{eqnarray}
Using the fact that under a trace over the space 1 one has
${\rm Tr}_1 \Big( R_{21} \cL_1 {R'}_{21} \Big) = {\rm Tr}_1
\Big( \cL_1 {R'_{21}}^{t_2} {R_{21}}^{t_2} \Big)^{t_2}$, one gets
\begin{equation}
  t(\beta_1) \, L_2(\beta_2)  
  =  
  L_2(\beta_2) \; {\rm Tr}_1
  \left[
  \cL_1(\beta_1) 
  \widehat R_{21}^*(\beta_2-\beta_1 + i\pi c)^{t_1t_2} 
  \left(\widehat R_{21}^* (\beta_2-\beta_1 + i\pi c )^{-1}
  \right)^{t_1t_2}
  \right]^{t_2}  \;.
\end{equation}
The last two terms in the right hand side cancel each other, leaving a trivial 
dependence in space 2 and ${\rm Tr}_1 \cL_1(\beta_1) \equiv t(\beta_1)$ in 
space 1.  This shows the commutation of $t(\beta_1)$ with $L(\beta_2)$, that 
is with the full algebra $\dey$ at $c=-N$.  
\finproof

\begin{coro}
  At the critical level $c=-N$ the operators $t(\beta)$ realise an
  Abelian algebra
  \begin{equation}
    t(\beta_1) \; t(\beta_2) = t(\beta_2) \; t(\beta_1) \;.
  \end{equation}
\end{coro}
The proof is obvious from the previous theorem. 

\subsection{Quantum exchange algebras}

\begin{lemm}
  When the parameters $r$ and $c$ satisfy the relation
  $Mr=-c-N$, with $M\in\ZZ$, the operators 
  $t(\beta) = {\rm Tr}[{\cal L}^{(M})(\beta)]$ realise an exchange algebra
  with the generators $L$ of $\dey$:
  \begin{equation}
    t(\beta_1) \; L(\beta_2) = 
    F(M,\beta_2-\beta_1) \; L(\beta_2) \; t(\beta_1) \;,
    \label{eq:qea}
  \end{equation}
  the function $F(M,\beta)$ being given by (\ref{eq:FM}). 
  \label{thm3}
\end{lemm}

\medskip
\noindent\textbf{Proof:}
The proof follows the lines of that of Theorem \ref{thm2}, 
starting from relation (\ref{eq31}) with $M\ne 0$. 
\finproof

\begin{thm}
  Under the condition $Mr=-c-N$, with $M\in\ZZ$, the operators 
  $t(\beta) = {\rm Tr}[{\cal L}^{(M)}(\beta)]$ close a quadratic algebra
  \begin{equation}
    t(\beta_1) \; t(\beta_2) = \cY_{N,r,M}(\beta_2-\beta_1) \;
    t(\beta_2) \; t(\beta_1) \;,
  \end{equation}
  where the function $\cY$ is given by 
  \begin{equation}
    {\cal Y}_{N,r,M}(\beta) = 
    \left\{ 
      \begin{array}{cc}
        \displaystyle \prod_{k=1}^{M} \quad
        \frac{\displaystyle \sin^2\frac{i\beta-k\pi r}{N}}
        {\displaystyle \sin^2\frac{i\beta+k\pi r}{N}}
        \;\;
        \frac{\displaystyle \sin\frac{i\beta+k\pi r+\pi}{N}}
        {\displaystyle \sin\frac{i\beta-k\pi r+\pi}{N}} 
        \;\;
        \frac{\displaystyle \sin\frac{i\beta+k\pi r-\pi}{N}}
        {\displaystyle \sin\frac{i\beta-k\pi r-\pi}{N}} 
        & \quad \mbox{for $M>0$} \,, \\ \\
        \displaystyle \prod_{k=1}^{|M|-1} \quad 
        \frac{\displaystyle \sin^2\frac{i\beta-k\pi r}{N}}
        {\displaystyle \sin^2\frac{i\beta+k\pi r}{N}}
        \;\;
        \frac{\displaystyle \sin\frac{i\beta+k\pi r+\pi}{N}}
        {\displaystyle \sin\frac{i\beta-k\pi r+\pi}{N}} 
        \;\;
        \frac{\displaystyle \sin\frac{i\beta+k\pi r-\pi}{N}}
        {\displaystyle \sin\frac{i\beta-k\pi r-\pi}{N}} 
        & \quad \mbox{for $M<0$} \,. \\
      \end{array} 
    \right.
    \label{eqY}
  \end{equation}
  \label{thm4}
\end{thm}

\medskip
\noindent\textbf{Proof:}
Lemma \ref{thm3} and definitions (\ref{eq235}) imply that 
\begin{subeqnarray}
  t(\beta_1) \; L^+(\beta_2) 
  &=& F(M,\beta_2-\beta_1-\half \, i\pi c) \;
  L^+(\beta_2) \; t(\beta_1)
  \\
  t(\beta_1) \; L^-(\beta_2)^{-1} 
  &=& F(M,\beta_2-\beta_1- i\pi r)^{-1} \;
  L^-(\beta_2)^{-1} \; t(\beta_1)  \;.
\end{subeqnarray}
Since 
$t(\beta)={\rm Tr}[{\cal L}^{(M)}(\beta)]=  
{\rm Tr}
\left[(h^{-M} L^+(\beta-\half\,i\pi c))^{t} \, 
  (L^-(\beta)^{-1})^t 
\right]$,
one gets 
\begin{equation}
  t(\beta_1) \; t(\beta_2) = \frac{F(M,\beta_2-\beta_1 - i\pi c) }
  {F(M,\beta_2-\beta_1- i\pi r)} \;
  t(\beta_2) \; t(\beta_1) \;.
\end{equation}
The function $\cY_{N,r,M}(\beta_2-\beta_1)$ being defined by the
above ratio of $F$ functions, its explicit expression (\ref{eqY}) 
is obtained from formula (\ref{eq:FM}).
\finproof

\medskip
\noindent
When $N=2$, the expression (\ref{eqY}) becomes 
\begin{equation}
  {\cal Y}_{2,r,M}(\beta) = 
  \left\{ 
    \begin{array}{cc}
      \displaystyle \prod_{k=1}^{M} \quad
      \frac{\displaystyle \tan^2\frac{i\beta-k\pi r}{2}}
      {\displaystyle \tan^2\frac{i\beta+k\pi r}{2}}
      & \quad \mbox{for $M>0$} \,, \\ \\
      \displaystyle \prod_{k=1}^{|M|-1} \quad 
      \frac{\displaystyle \tan^2\frac{i\beta-k\pi r}{2}}
      {\displaystyle \tan^2\frac{i\beta+k\pi r}{2}}
      & \quad \mbox{for $M<0$} \,. \\
    \end{array} 
  \right.
  \label{eq:Y2}
\end{equation}

\begin{coro}
  When both relations $Mr=-c-N$ and $2r=Nh$ hold, with
  $M,h\in\ZZ$, the function $\cY_{N,r,M}(\beta)$ is equal to
  one. The operators $t(\beta)$ therefore generate an Abelian algebra. 
\end{coro}

\medskip
\noindent\textbf{Proof:} Direct calculation. 
\finproof

\section{Poisson structures on $\dey$}
\setcounter{equation}{0}

The results of the previous section have shown that the operators $t(\beta)$ 
generate an Abelian subalgebra of $\dey$ when certain conditions on the 
parameters $c$, $r$ are fulfilled.  One can then naturally induce a 
Poisson structure on this Abelian algebra by considering the exchange 
relations between the $t(\beta)$ in the neighbourhood of the critical line 
$c=-N$ or the lines $Mr=-c-N$ eq.  (\ref{eq:surf}).

\subsection{Poisson structure on the centre at $c=-N$}

\begin{thm}\label{thm5}
The Poisson structure of the generators $t(\beta)$ around $c=-N$ reads:
\begin{equation}
  \bigg\{ t(\beta_{1}),t(\beta_{2}) \bigg\} =  
  \frac{\displaystyle 
    -\frac{2\pi}{N} \; 
    \sin^2\frac{\pi}{N}\;\cos\frac{i(\beta_{2}-\beta_{1})}{N}} 
  {\displaystyle \sin\frac{i(\beta_{2}-\beta_{1})}{N} 
    \;\sin\frac{i(\beta_{2}-\beta_{1})+\pi}{N} 
    \;\sin\frac{i(\beta_{2}-\beta_{1})-\pi}{N}} \quad
  t(\beta_{1}) \; t(\beta_{2}) 
  \,.
  \label{eqpoisbeta}
\end{equation}
\end{thm}

\medskip

\noindent \textbf{Proof:}

The proof follows the same scheme as the computation of the Poisson 
structures in the elliptic case at critical value $c=-N$.  We briefly
recall the basic steps: 

\medskip
\noindent
$\bullet$ Step 1: The exchange algebra for the generators $t(\beta)$
around  $c=-N$ reads
\begin{equation}
  t(\beta_{1}) \, t(\beta_{2}) = T(\beta_{2}-\beta_{1}) \; {\cal 
    M}(\beta_{2}-\beta_{1})^{i_1i_2}_{j_1j_2} \;\; {\cal
    L}(\beta_{2})^{j_2}_{i_2}  
  \;\; {\cal L}(\beta_{1})^{j_1}_{i_1} \,,
  \label{eq312}
\end{equation}
where
\begin{equation}
  {\cal M}(\beta) = \left(\left(R_{21}(\beta) \,
  {R_{21}(\beta-i\pi c-i\pi N)}^{-1} \, {R_{12}(-\beta)}^{-1}\right)^{t_2}
  \, {R_{12}(-\beta-i\pi c)}^{t_2}\right)^{t_2}
  \label{eq313}
\end{equation}
and
\begin{equation}
  T(\beta) = \frac{\displaystyle \sin\frac{i\beta-\pi}{N}}{\displaystyle 
    \sin\frac{i\beta+\pi}{N}} \; \frac{\displaystyle \sin\frac{i\beta+\pi-\pi 
      c}{N}}{\displaystyle \sin\frac{i\beta-\pi+\pi c}{N}} \;
    \frac{\displaystyle  
    \sin\frac{i\beta+\pi c}{N}}{\displaystyle \sin\frac{i\beta-\pi c}{N}} \;.
  \label{eq314}
\end{equation}
$\bullet$ Step 2: At $c=-N$, the derivative of ${\cal M}(\beta)$ with respect 
to $c$ vanishes due to the crossing-unitarity relation (it can actually be 
identified at $c=-N$ with the derivative of the crossing-unitarity
relation (\ref{eq:crossunit}) with  
respect to the spectral parameter).  Hence, the derivative of the exchange 
relation at $c=-N$ closes a Poisson structure on $t(\beta_{1})$.  \\
$\bullet$ Step 3: The Poisson structure function is now given by the 
derivative with respect to $c$ of the factor $T$ in eq.  (\ref{eq314}).  This 
factor is a product of the scalar normalisers used to rescale the structure 
$R$-matrix $\widehat R$ with respect to the unitary $R$-matrix $R$ 
in (\ref{eq226}).  Equation (\ref{eqpoisbeta}) immediately follows.
\finproof

\subsection{Poisson structure on the lines}

\begin{thm}
  Setting $Nh=2r+\epsilon$, 
  for any non zero integer $h$, one defines the $h$-labelled Poisson
  structure by 
  \begin{equation}
    \bigg\{ t(\beta_{1}),t(\beta_{2}) \bigg\}_h =  
    \lim_{\epsilon \rightarrow 0} \frac{1}{\epsilon} \;
    \Big( t(\beta_1)\; t(\beta_2) - t(\beta_2)\; t(\beta_1) \Big) \;.
    \label{eq:eps}
  \end{equation}
  Its explicit expression is 
  \begin{equation}
    \bigg\{ t(\beta_{1}),t(\beta_{2}) \bigg\}_h = f_h(\beta_2-\beta_1)
    \; t(\beta_{1}) \; t(\beta_{2})
  \end{equation}
  where
  \begin{equation}
    f_h(\beta) =  
    \cases{
    M(M+1) f_s(\beta) 
    & for $h$ even \cr\cr
    2E(\sfrac{M}{2}) (E(\sfrac{M}{2})+1) f_s(\beta)
    + 2E(\sfrac{M+1}{2})^2 f_c(\beta) 
    & for $h$ odd \cr
    }
    \label{eq:poish}
  \end{equation}
  with 
  \begin{equation}
    f_s(\beta) =  
    \frac{\displaystyle 
      -\frac{\pi}{N} \; 
      \sin^2\frac{\pi}{N}\;\cos\frac{i\beta}{N}} 
    {\displaystyle \sin\frac{i\beta}{N} 
      \;\sin\frac{i\beta+\pi}{N} 
      \;\sin\frac{i\beta-\pi}{N}} 
    \label{eq:fh}
  \end{equation}
  and
  \begin{equation}
    f_c(\beta) =  
    \frac{\displaystyle 
      \frac{\pi}{N} \; 
      \sin^2\frac{\pi}{N}\;\sin\frac{i\beta}{N}} 
    {\displaystyle \cos\frac{i\beta}{N} 
      \;\cos\frac{i\beta+\pi}{N} 
      \;\cos\frac{i\beta-\pi}{N}} 
    \label{eq:gh}
  \end{equation}
\end{thm}

\medskip
\noindent\textbf{Proof:}
By direct calculation, after noting that the right hand side of 
(\ref{eq:eps}) is equal to 
$\left. \displaystyle \frac{d \cY}{d \epsilon} \right|_{\epsilon=0}
t(\beta_1)t(\beta_2)$.

\finproof

\medskip
\noindent
In the case $N=2$, the functions $f_s$ and $f_c$ take the simple form
\begin{equation}
  f_s(\beta) = - f_c(\beta) = \frac{\pi}{\sin i\beta} \;.
  \label{eq:fg2}
\end{equation}
Therefore there is only one type of Poisson structure in the $sl(2)$
case since  the normalisation factors in the r.h.s. of (\ref{eq:poish})
may always be reabsorbed in the definition of $\epsilon$.

\section{Higher spin generators}

\begin{thm}\label{thmwl}
We define the operators $w_{s}(\beta)$ ($s = 1, \dots, N-1$) by:
\begin{equation}
w_{s}(\beta) \equiv \mbox{Tr} \left[ \prod_{1 \le k \le s}^{{\textstyle 
\curvearrowleft}} \left( \prod_{j>k} P_{kj} \right) \prod_{1 \le k \le 
s}^{{\textstyle \curvearrowleft}} \left( {\cal L}_{k}^{(M)}(\beta_k) 
\prod_{j>k} \widehat R_{kj}^{*}(\beta_k-\beta_j+i\pi N)^{t_kt_j} \right) 
\right] \;,
\label{eqH1}
\end{equation}
where 
\begin{equation}
{\cal L}_{k}^{(M)}(\beta) \equiv (h^{-M} \, L^+_{k}(\beta-\half i\pi 
c))^{t_k} \, \widetilde L^-_{k}(\beta) \equiv \underbrace{\II \otimes \dots 
\otimes \II}_{k-1} \otimes (h^{-M} \, L^+(\beta-\half i\pi c))^{t} \, 
\widetilde L^-(\beta) \otimes \underbrace{\II \otimes \dots \II}_{s-k}
\end{equation}
with $M \in \ZZ$, $\beta_k = \beta-i\pi(k-\frac{s+1}{2})$, and $P_{kj}$ is 
the permutation operator between the spaces $k$ and $j$ including the 
spectral parameters.  \\
On the lines $c=-N-Mr$, the operators $w_{s}(\beta)$ realise an exchange 
algebra with the generators $L(\beta')$ of $\dey$:
\begin{equation}
w_{s}(\beta) \, L(\beta') = F^{(s)}\Big(M,\beta'-\beta\Big) \, L(\beta') \, 
w_{s}(\beta) \,,
\label{eqH2}
\end{equation}
where
\begin{equation}
F^{(s)}\Big(M,\beta'-\beta\Big) = \prod_{k=1}^{s} 
F\Big(M,\beta'-\beta_k\Big) \,.
\end{equation}
\end{thm}

\medskip

\noindent \textbf{Proof:} The proof of Theorem \ref{thmwl} is based on the 
same algebraic arguments as in the case of the full elliptic algebra (see 
Theorem 6 of ref.  \cite{AFRS5}).  Namely, the exchange relation 
(\ref{eqH2}) thus follows from the basic relation (\ref{eq31}), the
Yang--Baxter  
equation and the crossing-symmetry property (\ref{eq223}).    
\finproof

\medskip

\noindent One obtains the following corollary, the proof being obvious:
\begin{coro}
On the lines $c=-N-Mr$, the operators $w_{s}(\beta)$ realise an exchange 
algebra
\begin{equation}
w_s(\beta) \, w_{s'}(\beta') = \prod_{u=-\frac{s-1}{2}}^{\frac{s-1}{2}}
\prod_{v=-\frac{s'-1}{2}}^{\frac{s'-1}{2}}{\cal Y}_{N,r,M} 
\left(\beta'-\beta+i\pi(u-v)\right) w_{s'}(\beta') \, w_s(\beta) \,.
\end{equation}
When a additional relation $2r=Nh$ with $h \in \ZZ \backslash \{0\}$ is 
imposed, the function ${\cal Y}_{N,r,M}$ is equal to 1: the operators 
$w_{s}(\beta)$ realise then an Abelian subalgebra in $\dey$.
\end{coro}

\medskip

The previous corollary allows us to define Poisson structures on the 
corresponding Abelian subalgebras.  As usual, they are obtained as limits of 
the exchange algebras when $Nh=2r$ with $h \in \ZZ \backslash \{0\}$:
\begin{thm}\label{thmclass}
Setting $Nh=2r+\epsilon$ for any integer $h \ne 0$, the $h$-labelled Poisson 
structure defined by:
\begin{equation}
\{ w_{s}(\beta) , w_{s'}(\beta') \}_{(h)} = \lim_{\epsilon \rightarrow 0} 
\frac{1}{\epsilon} \, \Big( w_{s}(\beta) \, w_{s'}(\beta') - w_{s'}(\beta') \, 
w_{s}(\beta) \Big)
\end{equation}
has the following expression:
\begin{equation}
\{ w_{s}(\beta),w_{s'}(\beta') \} = \sum_{u=-(s-1)/2}^{(s-1)/2} 
\sum_{v=-(s'-1)/2}^{(s'-1)/2} f_h \Big(\beta'-\beta+i\pi(u-v)\Big) \, 
w_{s}(\beta) \, w_{s'}(\beta') \;,
\end{equation}
where
$f_h(\beta)$ is given by (\ref{eq:poish}).
\end{thm}

\medskip

\noindent \textbf{Proof:} The proof is made by direct calculation.
\finproof

\section{Conclusion}
We have proved that the existence of exchange subalgebras, Abelian
subalgebras and the Poisson structures on these objects, all survived
in the scaling limit $q\rightarrow 1$, 
$\frac{\ln p}{\ln q}\rightarrow r$,
$\frac{\ln z}{\ln q}\rightarrow \frac{i\beta}{\pi}$,
and the characteristic relations and structure functions were given by
the scaling limits, suitably defined whenever potential divergences
arise, of their analogues in the full elliptic case.\\
Two directions should now be investigated. We have already commented
upon the necessary careful study of precise mode expansions of the
generating functionals and structure function, required to construct
explicit representations of such algebras. \\
In a more abstract setting, the question arises of interpreting more
precisely the ``deformation'' involved when going from $\deysr$ (at
$r=\infty$) to $\dey$. We have established at $N=2$ that this
deformation is in fact a Drinfel'd twist generated by the
$q\rightarrow 1$ scaling limit of the twist
${\cU}_{q}(\widehat{sl}(2)_{c})\rightarrow
{\cA}_{q,p}(\widehat{sl}(2)_{c})$  \cite{JKOS}. We shall report on
this and other very interesting connections established between
various algebraic objects of similar type \cite{AAFRRS}.

\bigskip

\noindent \textbf{Acknowledgements}

This work was supported in part by CNRS and EC network contract number 
FMRX-CT96-0012.  M.R. was supported by an EPSRC research grant 
no. GR/K 79437. 
The work greatly benefitted from penetrating observations and
suggestions of P.~Sorba. 
The authors would like to thank V.~Korepin for his stimulating
questions and S.M.~Khoroshkin and H.~Konno for
valuable clarifications. 
J.A. and M.R. wish to thank the LAPTH for its kind 
hospitality.

\end{document}